\documentclass[a4paper]{amsart}

\usepackage[margin=1.2in]{geometry}
\usepackage{bm}
\usepackage{array}
\usepackage{float}
\usepackage{listings}
\usepackage{xcolor}
\usepackage{tabularray}

\newtheorem{theorem}{Theorem}
\theoremstyle{definition}
\newtheorem{definition}[theorem]{Definition}
\newtheorem{remark}[theorem]{Remark}
\newtheorem{lemma}[theorem]{Lemma}
\newtheorem{examplex}[theorem]{Example}
\newenvironment{example}
  {\pushQED{\qed}\examplex}
  {\popQED\endexamplex}

\definecolor{codegreen}{rgb}{0,0.6,0}
\definecolor{codegray}{rgb}{0.5,0.5,0.5}
\definecolor{sageblue}{RGB}{94, 94, 255}
\definecolor{codepurple}{rgb}{0.58,0,0.82}
\definecolor{backcolour}{rgb}{0.95,0.95,0.92}

\lstdefinestyle{mystyle}{
    frame=tb,
    keywordstyle=\color{magenta},
    escapeinside={@}{@},
    backgroundcolor=\color{backcolour},   
    commentstyle=\color{codegreen},
    numberstyle=\tiny\color{codegray},
    stringstyle=\color{codegreen},
    basicstyle=\ttfamily\footnotesize,
    breakatwhitespace=false,         
    breaklines=true,                 
    captionpos=b,                    
    keepspaces=true,                 
    numbersep=5pt,                  
    showspaces=false,                
    showstringspaces=false,
    showtabs=false,                  
    tabsize=2,
}

\def\one{\mathbf{1}} 
\def\zero{\mathbf{0}} 
\bmdefine{\ii}{i}
\bmdefine{\rr}{r} 
\bmdefine{\ww}{w}
\bmdefine{\yy}{y}
\bmdefine{\zz}{z}

\newcommand{\Z}{{\mathbb Z}}
\newcommand{\Q}{{\mathbb Q}}
\newcommand{\C}{{\mathbb C}}
\newcommand{\R}{{\mathbb R}}
\newcommand{\N}{{\mathbb N}}
\newcommand{\dom}{\mathcal{D}}
\newcommand{\sing}{\mathcal{V}}
\newcommand{\hes}{\mathcal{H}}
\newcommand{\code}[1]{\texttt{\detokenize{#1}}}

\usepackage[backend=bibtex]{biblatex}
\begin{document}

\title[Rigorous ACSV in SageMath]{Rigorous Analytic Combinatorics in Several Variables\\ in SageMath}

\author[B. Hackl]{Benjamin Hackl}
  \address{
    Department of Mathematics, University of Klagenfurt, Austria}
  \curraddr{
    Institute of Mathematics and Scientific Computing, University of Graz, Austria}
  \email{benjamin.hackl@uni-graz.at}

\author[A. Luo]{Andrew Luo}
  \address{Department of Combinatorics and Optimization, University of Waterloo, Canada}
  \email{j92luo@uwaterloo.ca}

\author[S. Melczer]{Stephen Melczer}
  \address{Department of Combinatorics and Optimization, University of Waterloo, Canada}
  \email{smelczer@uwaterloo.ca}

\author[J. Selover]{Jesse Selover}
  \address{Department of Mathematics, University of Massachusetts, Amherst, USA}
  \email{jselover@umass.edu}

\author[E. Wong]{Elaine Wong}
  \address{Oak Ridge National Laboratory, Tennessee, USA}
  \email{wongey@ornl.gov}
  
\begin{abstract}
We introduce the new \textbf{sage\_acsv} package for the SageMath computer algebra system, allowing users to rigorously compute asymptotics for a large variety of multivariate sequences with rational generating functions. Using Sage's support for exact computations over the algebraic number field, this package provides the first rigorous implementation of algorithms from the theory of analytic combinatorics in several variables.
\end{abstract}

\maketitle

\section{Introduction}

The field of \emph{Analytic Combinatorics in Several Variables (ACSV)}~\cite{Melczer2021,PemantleWilson2013} adapts methods from complex analysis in several variables, algebraic and differential geometry, topology, and computer algebra to create effective methods to study the asymptotic properties of a multivariate sequence $(f_{\ii})_{\ii \in \N^d}$ using analytic properties of its \emph{generating function}, 
\[ F(\zz) = \sum_{\ii \in \N^d}f_{\ii}\zz^{\ii} = \sum_{(i_1,\ldots,i_d)\in\N^d}f_{i_1,\ldots, i_d} \; z_1^{i_1}\cdots z_d^{i_d} , \]
when this series represents a complex analytic function near the origin. 
The \textbf{sage\_acsv} package\footnote{A demonstration of this work was given at the 35th International Conference on Formal Power Series and Algebraic Combinatorics in Davis, CA, USA and is published in the corresponding proceedings. The published version of this paper can be found here \url{https://www.mat.univie.ac.at/~slc/wpapers/FPSAC2023/90.pdf}.} takes a multivariate rational function $F(\zz)$ specified as an explicit symbolic fraction
and a \emph{direction vector} $\rr \in \Z_{>0}^d$ and returns an asymptotic expansion for the \emph{$\rr$-diagonal} sequence\footnote{The theory of ACSV shows how asymptotics typically vary uniformly with small perturbations of $\rr$, so it is also possible to derive asymptotics varying near fixed directions, and even limit theorems, using similar techniques.} of the power series coefficients $(f_{n\rr})_{n \geq 0}$ of $F$, whose indices are natural number multiples of $\rr$. Our package works under conditions that are verifiable and hold generically, except that $F$ must be \emph{combinatorial}, meaning that all but a finite number of its power series coefficients are non-negative. 

\begin{remark}
It is unknown whether combinatorality is a computationally decidable property, even in the univariate case. However, this property holds often in combinatorial contexts (for instance, the multivariate generating function of any combinatorial class is combinatorial because its series coefficients count something).
\end{remark}

The \code{sage_acsv} package can be installed in any recent SageMath installation
(preferably version \texttt{9.4} or later) by running the command \code{sage -pip install sage-acsv}
from any Sage instance with access to the internet, or by downloading its source code at \url{https://github.com/ACSVMath/sage_acsv} and placing the module in the appropriate Python search path.

\begin{example}
\label{ex:1}
The $(1,1)$-diagonal of the combinatorial rational function
\[F(x,y) = \frac{1}{1-x-y} = \sum_{a,b \geq 0}\binom{a+b}{a}x^ay^b\] 
forms the sequence $f_{n,n} = \binom{2n}{n}$ of central binomial coefficients. After installing the package, running the code

\begin{lstlisting}[language=python, showstringspaces=false]
@\color{sageblue}sage:@ from sage_acsv import diagonal_asy
@\color{sageblue}sage:@ var('x, y')
(x, y)
@\color{sageblue}sage:@ diagonal_asy(1/(1 - x - y), as_symbolic=True)
4^n/(sqrt(pi)*sqrt(n))
\end{lstlisting}

\noindent
verifies that the required assumptions of ACSV hold and proves that
\[ 
\binom{2n}{n} 
= \frac{4^n}{\sqrt{\pi n}}\left(1 + O\left(\frac{1}{n}\right)\right).
\]
\end{example}

The asymptotic behavior of a multivariate rational diagonal sequence under our assumptions (to be detailed below) is specified by a sum of terms of the form $C \rho^n (\pi n)^\alpha$ where $C$ and $\rho$ are algebraic numbers and $\alpha$ is a rational number. By default, the \code{diagonal_asy} command returns expansions of this form as a list of tuples containing elements of Sage's \code{Algebraic Field}, which stores exact representations of the quantities. The optional flag \code{as_symbolic} invoked in Example~\ref{ex:1} tells \code{diagonal_asy} to return a symbolic sum involving $n$ which makes it easier to view (and approximate) asymptotics but makes it difficult to access the exact algebraic numbers involved.

\begin{example}
\label{ex:2}
The asymptotic behavior of the sequence
\[ f_n = \sum_{k=0}^n\binom{n}{k}^2\binom{n+k}{k}^2, \] 
which is the main diagonal of 
\[F(w,x,y,z) = \frac{1}{1 - w(1+x)(1+y)(1+z)(xyz+yz+y+z+1)}\ ,\] 
was used by Apéry~\cite{Apery4} in his celebrated proof of the irrationality of $\zeta(3)$. Running

\begin{lstlisting}[language=python, showstringspaces=false]
@\color{sageblue}sage:@ from sage_acsv import diagonal_asy
@\color{sageblue}sage:@ var('w, x, y, z, t')
(w, x, y, z, t)
@\color{sageblue}sage:@ G = 1
@\color{sageblue}sage:@ H = 1 - w*(1 + x)*(1 + y)*(1 + z)*(1 + y + z + y*z + x*y*z)
@\color{sageblue}sage:@ diagonal_asy(G/H, as_symbolic=True)
1.225275868941647?*33.97056274847714?^n/(pi^1.5*n^1.5)
\end{lstlisting}
\noindent
gives the dominant asymptotic behavior of $f_n$ in terms of algebraic numbers represented by decimal approximations. Running

\begin{lstlisting}[language=python, showstringspaces=false]
@\color{sageblue}sage:@ asm_vals = diagonal_asy(G/H)
@\color{sageblue}sage:@ add([a.radical_expression()^n*b*c*d.radical_expression()
          for (a, b, c, d) in asm_vals])
(12*sqrt(2) + 17)^n*sqrt(17/2*sqrt(2) + 12)/(4*pi^1.5*n^1.5)
\end{lstlisting}
\noindent
represents these algebraic quantities in radicals using Sage's capability
for computing symbolic radical expressions, proving
\[ f_n = \frac{\left(12\sqrt{2}+17\right)^n}{\sqrt{(\pi n)^3}}\left(\frac{\sqrt{12+17/\sqrt{2}}}{4} + O\left(\frac{1}{n}\right)\right). \]
We note that such radical expressions are not certified and, of course, it is not always possible (or, even if possible, useful) to write the algebraic numbers appearing in asymptotics for general multivariate diagonals in radicals. 
\end{example}

The, by now well-established, theory of univariate analytic combinatorics~\cite{FlajoletSedgewick2009} shows that to determine asymptotics of the power series coefficients of a generating function $f(z)$, one should find the \emph{dominant singularities} of $f$ (those closest to the origin), determine the type of singularity (pole of finite order, logarithmic or algebraic branch point, essential singularity, etc.), and then use known \emph{transfer theorems} to determine the contribution of each singularity to the asymptotics for the corresponding sequence. Dominant singularities are generalized from one to several variables by the following definition.

\begin{definition}[minimal point]
Let $F(\zz) = \sum_{\ii\in\N^d}f_\ii\zz^\ii$ be a convergent power series with domain of convergence $\dom$. A singularity $\ww\in\R^d$ of $F(\zz)$ is said to be \emph{minimal} if $\ww$ lies in the boundary $\partial \dom$. Equivalently, $\ww$ is minimal if there does not exist another singularity $\ww' \in \R^d$ of $F(\bm{z})$ with $|w_i'| < |w_i|$ for all $1 \le i \le d$.
\end{definition}

Although the theory of ACSV applies to meromorphic functions, we focus here on the case when $F(\zz)=G(\zz)/H(\zz)$ is a rational function so that we can use tools from polynomial system solving (rational functions can also capture the behavior of algebraic functions through embeddings in higher dimensions~\cite{GreenwoodMelczerRuzaWilson2022}). We always assume that $G$ and $H$ are coprime polynomials, which implies that the set of singularities of $F$ is defined by the \emph{singular variety} $\sing = \{\zz\in\C^d : H(\zz)=0\}$.

Unlike the univariate case, in at least two variables there are always an infinite number of minimal points when $F$ admits singularities. This reflects, in part, the fact that minimal points are important to the asymptotics of $(f_\ii)$ and there are an infinite number of ways for the index $\ii$ to approach infinity by varying the direction vector $\rr$. Thankfully, it is possible in generic circumstances to determine a finite set of singularities, depending on $\rr$, that dictate asymptotics of the $\rr$-diagonal of $F$.

\begin{definition}[critical point]
A point $\ww \in \C_*^d := (\C\setminus \{0\})^d$ with nonzero coordinates is a \emph{(simple smooth) critical point} of $F$ for the direction $\rr=(r_1,\dots,r_d)$ if the gradient $(\nabla H)(\ww)$ is nonzero and
\begin{align}
 \begin{split}
    H(\ww)&=0,\\
    r_kz_1H_{z_1}(\ww) - r_1z_kH_{z_k}(\ww) &= 0 \; \text{ for } 2\leq k\leq d,
    \end{split}  \label{eq:CP}
\end{align}
where $H_{z_k}$ denotes the partial derivative of $H$ with respect to the variable $z_k$.
\end{definition}

The theory of ACSV shows how minimal critical points, when they exist, typically determine asymptotics (see Theorem~\ref{thm:ACSV} below). The command 

\begin{center}
	\code{diagonal_asy(F, r, linear_form, return_points, as_symbolic)}
\end{center}

\noindent
computes asymptotics of the $\rr$-diagonal $f_{n\rr}$ of $F(\zz)=G(\zz)/H(\zz)$ as $n\rightarrow\infty$, under the following assumptions:

\begin{enumerate}
    \item $H(\zero)\neq0$, so that $G(\zz)/H(\zz)$ has a convergent power series expansion near the origin,
    \item $H$ and all of its partial derivatives do not simultaneously vanish, so that all poles of $F$ are simple and $\sing$ forms a manifold,
    \item the smooth critical point system~\eqref{eq:CP} has a finite number of solutions, at least one of which is minimal, and
    \item the explicit matrix $\hes=\hes(\ww)$ defined by~\eqref{eq:hes} below is nonsingular at all minimal critical points $\ww$.
\end{enumerate}
The algorithm verifies these assumptions computationally, so the user need not worry about them. As mentioned above, we further require that
\begin{enumerate}
\item[5.] the power series expansion of $F$ is \emph{combinatorial}, meaning that it admits at most a finite number of negative coefficients,
\end{enumerate}
however it is unknown whether combinatorality is decidable, even for univariate rational functions, so \textbf{the user must know this} through other means in order for the output to be proven correct (for instance, if the input is a multivariate generating function for a counting sequence then this is satisfied). Table \ref{tab:parameters} summarizes the possible input parameters of the function \code{diagonal_asy}.

Multivariate generating functions enumerating combinatorial classes are, by their definition, combinatorial. The form of $F$ can also be used to prove combinatorality, for instance if $F(\zz)=G(\zz)/(1-I(\zz))$ where $G$ and $I$ are polynomials with non-negative coefficients and $I$ vanishes at the origin.

\begin{remark}
    If $H$ and its partial derivatives simultaneously vanish because $H$ is not square-free then minimal critical points can still be determined by replacing $H$ with its square-free part (the product of its irreducible factors); asymptotics can also be determined using a minor generalization of the formula given below to higher order poles. If $\sing$ admits points where it is not locally a manifold, but~\eqref{eq:CP} still has a finite set of solutions and admits a smooth minimal critical point $\ww$ such that all other solutions of~\eqref{eq:CP} with the same coordinate-wise modulus are also smooth, then this approach can still determine asymptotics. In general, if nonsmooth minimal points affect asymptotic behavior, then more advanced results are required~\cite[Part III]{Melczer2021}.
\end{remark}

\begin{table}[ht]
\centering
\caption{Parameters for the \code{diagonal_asy} command.}
\label{tab:parameters}
\begin{tblr}{ |p{0.18\linewidth}< {\centering} |p{0.28\linewidth}< {\centering} |p{0.38\linewidth} < {\centering}|}
 \hline 
Parameter & Type & Description \\ \hline \hline
$F(\zz)$ & Combinatorial rational function in $\Z(z_1, ..., z_d)$ & Ratio of coprime polynomials $G(\zz)$ and $H(\zz)$.\\ \hline
$\rr$  (optional) & List of $d$ positive integers & Direction to compute asymptotics. If none is provided, use $\rr=\one$. \\ \hline
\code{linear_form} (optional) & Linear polynomial in $\Z[z_1, \dots, z_d]$ & Integer linear form to be used in the algorithm. If none provided, generate one at random. \\  \hline
\code{return_points} (optional) & True or False & Flag to also return coordinates of minimal critical points determining asymptotics. Default is False. \\  \hline
\code{as_symbolic} (optional) & True or False & Flag to return asymptotics as a symbolic sum involving $n$. Default is False. \\ 
 \hline
\end{tblr}
\end{table}

\begin{example}
The sequence alignment problem in molecular biology compares evolutionary relationships between species by measuring differences in their DNA sequences. Pemantle and Wilson~\cite{PemantleWilson2008} study one such problem whose behavior is encoded by the main diagonal of the combinatorial multivariate rational generating function
\[ \frac{x^2y^2-xy+1}{1-(x+y+xy-xy^2-x^2y+x^2y^3+x^3y^2)}.\]
We immediately prove dominant asymptotic behavior for the main diagonal by running

\begin{lstlisting}[language=python, showstringspaces=false]
@\color{sageblue}sage:@ from sage_acsv import diagonal_asy
@\color{sageblue}sage:@ var('x, y')
(x, y)
@\color{sageblue}sage:@ G = x^2*y^2 - x*y + 1
@\color{sageblue}sage:@ H = 1 - (x + y + x*y - x*y^2 - x^2*y + x^2*y^3 + x^3*y^2)
@\color{sageblue}sage:@ F = G/H
@\color{sageblue}sage:@ asm_vals = diagonal_asy(F, as_symbolic=True)
0.9430514023983397?*4.518911369262258?^n/(sqrt(pi)*sqrt(n))
\end{lstlisting}

\noindent
where now the leading constant is an algebraic number of degree 10 and the exponential growth constant is an algebraic number of degree 5.
\end{example}

A test notebook working through these examples, and further ones, is available with the source code
of the package at \url{https://github.com/ACSVMath/sage_acsv}.

\subsubsection*{Past Work}

The theoretical underpinnings of our software began with the work of Pemantle and Wilson~\cite{PemantleWilson2002} on smooth ACSV. This early paper used easy to understand explicit contour deformations of complex integrals to derive asymptotics, however it relies on stronger assumptions than we need. The specific asymptotic result we use was originally proven in Baryshnikov and Pemantle~\cite{BaryshnikovPemantle2011} using the notion of \emph{hyperbolic cones} to enable more advanced deformations, which was then studied from an algorithmic viewpoint by Melczer and Salvy~\cite{MelczerSalvy2021}. The latter two authors created a preliminary Maple implementation of the algorithm which was not rigorous, as it did not use certified numerics. 
A previous implementation by Raichev~\cite{Raichev2011}, currently included as a core module in SageMath,
can also compute\footnote{Unfortunately, as of the publication of this article, changes to the underlying SageMath codebase have broken some functionality of Raichev's package.} the asymptotic contributions of minimal critical points (including in some nonsmooth situations). However, this package requires the user to independently find and certify the minimal critical points (the hardest step of the analysis), so it also does not rigorously determine asymptotics. In contrast, our package uses Sage's exact computations for algebraic numbers to rigorously decide the (in)equalities necessary to prove computed asymptotics under our assumptions.

\section{Algorithmic Details}

In this section we quickly recap the theoretical background of ACSV, and then describe how our package works in more detail.

\subsection{Recap of ACSV}

We start by representing the $\rr$-diagonal sequence by a Cauchy integral,
\[ f_{n\rr} = \frac{1}{(2\pi i)^d} \int_{\mathcal{C}} F(\zz) \frac{d\zz}{\zz^{n\rr+\one}}, \]
where $\mathcal{C}$ is any product of circles sufficiently close to the origin, with $d\zz=dz_1\cdots dz_d$ and $\zz^{n\rr+\one}=(z_1^{nr_1+1},\ldots, z_d^{nr_d+1})$. Under our assumptions, it is possible to deform $\mathcal{C}$ away from the origin and replace the Cauchy integral by a residue integral localized to smooth minimal critical points (except at points that yield an exponentially negligible error). The definition of critical points is crafted so that the resulting residue integral can be analyzed using the classical saddle-point method. Our asymptotic formulas depend on a certain matrix, which we now define.

\begin{definition}[phase Hessian matrix] 
    If $\ww$ is a smooth critical point, then the \emph{phase Hessian matrix} $\hes=\hes(\ww)$ at $\ww$ is the $(d-1)\times(d-1)$ matrix defined by
    \begin{equation}
    \hes_{i,j} = 
    \begin{cases}
        V_iV_j + U_{i,j} - V_jU_{i,d} - V_iU_{j,d} + V_iV_jU_{d,d}\, , & i \neq j, \\[+3mm]
        V_i + V_i^2 + U_{i,i} - 2V_iU_{i,d} + V_i^2U_{d,d}\, , & i=j,
    \end{cases}
    \label{eq:hes}
    \end{equation}
    where
    \[ U_{i,j} = \frac{w_iw_j H_{z_iz_j}(\ww)}{w_dH_{z_d}(\ww)} \qquad \text{and} \qquad V_i = \frac{r_i}{r_d}.\]
\end{definition}

The key ACSV theorem in our context is the following.

\begin{theorem}[{Melczer~\cite[Theorem 5.1]{Melczer2021}}]
    \label{thm:ACSV}
    Suppose that the system of polynomial equations~\eqref{eq:CP} admits a finite number of solutions, exactly one of which, $\ww\in\overline{\Q}_*^d$, is minimal. Suppose further that $H_{z_d}(\ww)\neq0$, and that the phase Hessian matrix $\hes$ at $\ww$ has nonzero determinant. Then, as $n\rightarrow\infty$,
    
    {\small
        \[
        f_{n\rr} = \ww^{-n\rr} n^{(1-d)/2} \frac{(2\pi r_d)^{(1-d)/2}}{\sqrt{\det(\hes)}} \left(\frac{-G(\ww)}{w_d\, H_{z_d}(\ww)} + O\left(\frac{1}{n}\right)\right).
        \]
    }
    
    If the zero set of $H$ contains a finite number of points with the same coordinate-wise modulus as $\ww$, all of which satisfy the same conditions as $\ww$, then an asymptotic expansion of $f_{n\rr}$ is obtained by summing the right hand side of this expansion at each point.
\end{theorem}

\begin{remark}
    If $G(\ww)=0$, then the leading asymptotic term in Theorem~\ref{thm:ACSV} will vanish. In many cases dominant asymptotics can usually still be determined by computing higher-order terms using (increasingly complicated) explicit formulas, however (in nongeneric situations, or when $G$ and $H$ were not reduced to be coprime) it is possible for $f_{n\rr}$ to grow exponentially slower than~$\ww^{-n\rr}$, and a local analysis near $\ww$ can only prove $f_{n\rr} = O(\ww^{-n\rr}n^{-k})$ for any positive integer $k$.
\end{remark}

In order to determine which critical points are minimal we use the following result.

\begin{lemma}[Melczer and Salvy~\cite{MelczerSalvy2021}]
    \label{lem:combmin}
    If $F$ is combinatorial and $\yy\in\C_*^d$ is a minimal critical point then so is $(|y_1|,\dots,|y_d|)$. Furthermore, $\ww \in \R_{>0}^d$ is a minimal critical point if and only if the system
    \begin{equation}
        \begin{split}
            H(\zz) = H(tz_1,\dots,tz_d) &= 0, \\
            z_1H_{z_1}(\zz)-r_1\lambda = \cdots = z_dH_{z_d}(\zz)-r_d\lambda &= 0,
        \end{split}
        \label{eq:extendedSys}
    \end{equation}
    has a solution $(\zz,\lambda,t) \in \R^{d+2}$ with $\zz=\ww$ and $t=1$ \emph{and} no solution with $\zz=\ww$ and $0 < t < 1$.
\end{lemma}

Theorem~\ref{thm:ACSV} and Lemma~\ref{lem:combmin} display our overall strategy: encode the (generically finite~\cite[Section~5.3.4]{Melczer2021}) set of solutions to the polynomial system~\eqref{eq:extendedSys} in a convenient manner, determine which solutions with $t=1$ have $z_1,\dots,z_d$ positive and real, and then use the solutions with $t \in (0,1)$ to eliminate those that are not minimal. Under our assumptions there will be at most one positive real minimal critical point: it remains to find the other critical points with the same coordinate-wise modulus and then add the asymptotic contributions of each given by Theorem~\ref{thm:ACSV}.

\subsection{Kronecker Representations}

The key to an efficient algorithm is the representation used to encode the solutions of the extended critical point system~\eqref{eq:extendedSys}. Following Melczer and Salvy~\cite{MelczerSalvy2021}, we use a \emph{Kronecker representation}, which is also known as a \emph{rational univariate representation} (see Melczer and Salvy~\cite{MelczerSalvy2021} for background on the Kronecker representation and its history in computer algebra). A \emph{Kronecker representation} of a zero-dimensional algebraic set 
\[ S = \{\zz \in \mathbb{C}^d : f_1(\zz) = ... = f_d(\zz) = 0\}\] 
consists of
\begin{itemize}
    \item a new variable $u$ given by a separating integer linear form 
    $u = \boldsymbol{\kappa} \cdot \zz$ in the original variables $\zz$ for some $\boldsymbol{\kappa}\in \mathbb{Z}^d$,
    \item a square-free integer polynomial $P \in \Z[u]$, and
    \item integer polynomials $Q_1, ..., Q_d \in \Z[u]$,
\end{itemize}
such that the points $\zz$ defined by $z_i=Q_i(u)/P'(u)$ as $u$ runs through the roots of $P$ give the elements of $S$. The Kronecker representation of a zero-dimensional variety encodes its points using the univariate polynomial $P$, and is constructed so that the degrees and maximum coefficient sizes of $P$ and the $Q_i$ can be efficiently bounded~\cite{Safey-El-DinSchost2018,Schost2001}.

Although specialized, Gröbner-free, algorithms to compute a Kronecker representation~\cite{GiustiLecerfSalvy2001} exist, to the best of our knowledge they have not been implemented in Sage. We thus use a lexicographical Gröbner basis computation to determine our Kronecker representations. The ability to compute Kronecker representations in Sage may be of independent interest to some users.

\begin{example}
Determining the asymptotics of the Apéry-3 Sequence requires solving an \emph{extended critical point system}
\begin{align*}
    xyz - xz + x - \lambda &= 0, \\
    xyz - yz + y - \lambda &= 0, \\
    xyz - xz - yz + z - \lambda &= 0, \\
    xyz - xz - yz + x + y + z - 1 &= 0, \\
    xyzt^3 - xzt^2 - yzt^2 + xt + yt + zt - 1 &= 0.
\end{align*}
Using the linear form $u = x+t$, running
\begin{lstlisting}[language=python, showstringspaces=false]
@\color{sageblue}sage:@ from sage_acsv import kronecker
@\color{sageblue}sage:@ var('x, y, z, t, lambda_')
(x, y, z, t, lambda_)
@\color{sageblue}sage:@ kronecker([x*y*z - x*z + x - lambda_, x*y*z - y*z + y - lambda_,
....:           x*y*z - x*z - y*z + z - lambda_,
....:           x*y*z - x*z - y*z + x + y + z - 1,
....:           x*y*z*t^3 - x*z*t^2 - y*z*t^2 + x*t + y*t + z*t -1],
....:          [x, y, z, t, lambda_], x + t)
(u_^8 - 18*u_^7 + 146*u_^6 - 692*u_^5 + 2067*u_^4 - 3922*u_^3 + 4553*u_^2 - 2925*u_ + 790,
 [10*u_^7 - 153*u_^6 + 1046*u_^5 - 4081*u_^4 + 9589*u_^3 - 13270*u_^2 + 9844*u_ - 2985,
  10*u_^7 - 154*u_^6 + 1061*u_^5 - 4180*u_^4 + 9954*u_^3 - 14044*u_^2 + 10714*u_ - 3380,
  -u_^7 + 11*u_^6 - 56*u_^5 + 157*u_^4 - 182*u_^3 - 140*u_^2 + 527*u_ - 335,
  8*u_^7 - 139*u_^6 + 1030*u_^5 - 4187*u_^4 + 10021*u_^3 - 14048*u_^2 + 10631*u_ - 3335,
  -12*u_^7 + 181*u_^6 - 1231*u_^5 + 4801*u_^4 - 11275*u_^3 + 15548*u_^2 - 11452*u_ + 3440])
\end{lstlisting}
gives the Kronecker representation
\begin{align*}
    P(u) =& u^8 - 18u^7 + 146u^6 - 692u^5 + 2067u^4 - 3922u^3 + 4553u^2 - 2925u + 790, \\
    Q_x(u) =& 10u^7 - 153u^6 + 1046u^5 - 4081u^4 + 9589u^3 - 13270u^2 + 9844u - 2985, \\
    Q_y(u) =& 10u^7 - 154u^6 + 1061u^5 - 4180u^4 + 9954u^3 - 14044u^2 + 10714u - 3380, \\
    Q_z(u) =& -u^7 + 11u^6 - 56u^5 + 157u^4 - 182u^3 - 140u^2 + 527u - 335, \\
    Q_t(u) =& 8u^7 - 139u^6 + 1030u^5 - 4187u^4 + 10021u^3 - 14048u^2 + 10631u - 3335, \\
    Q_\lambda(u) =& -12u^7 + 181u^6 - 1231u^5 + 4801u^4 - 11275u^3 + 15548u^2 - 11452u + 3440.
    \qedhere
\end{align*}
\end{example}

After computing a Kronecker representation of the critical point system, we use Sage's solver over the \code{Real Algebraic Field} to determine the real roots of $P(u)$, identify which correspond to critical points with positive real coordinates, filter out those that are not minimal by examining solutions with $t$ coordinate in $(0,1)$, and then identify which critical points have the same coordinate-wise modulus. All identities and inequalities are verified exactly by working over the \code{Algebraic Field} in Sage (numeric approximations to sufficient accuracy are used to decide inequalities, unless symbolic computations are absolutely necessary to prove equality). Further details can be found in the documentation for our package.

\subsection{Future Work}

There are many ways to generalize this work, including: computing the higher-order terms in the expansion given by Theorem~\ref{thm:ACSV}, working with nonsmooth critical points, Gröbner-free methods of computing Kronecker representations, and generalizing to larger classes of meromorphic functions where computations are still effective. From a complexity point of view, the most expensive operation to find asymptotics is typically identifying critical points with the same coordinate-wise modulus. Our package currently uses a naive method relying on Sage's built-in support for algebraic numbers, however work is ongoing to adapt a more efficient method described in Melczer and Salvy~\cite{MelczerSalvy2021} in the specialized ACSV setting. Furthermore, we use Sage’s interface to Singular for necessary Gröbner basis computations. Work adapting newer packages, such as \code{msolve}~\cite{Berthomieu2021}, is also ongoing.

\section*{Acknowledgements}
The authors acknowledge the support of the AMS Math Research Community \emph{Combinatorial Applications of Computational Geometry and Algebraic Topology}, which was funded by the National Science Foundation under Grant Number DMS 1641020. SM and AL's work partially supported by NSERC Discovery Grant RGPIN-2021-02382 and the NSERC Undergraduate Summer Research Program.

\printbibliography

\end{document}